\documentclass{amsart}

\usepackage{graphicx}

\def\co{\colon\thinspace}
\newcommand{\R}{\mathbb{R}}
\newcommand{\D}{\mathcal{D}}
\newcommand{\ogamma}{\overline{\gamma}}
\newcommand{\ogF}{\ogamma_{\mathrm{F}}}
\newcommand{\rot}{{\tt rot}}

\theoremstyle{definition}
\newtheorem*{thm}{Theorem}

\begin{document}

\title{Horizontal Loops in Engel Space}

\author{Hansj\"org Geiges}

\address{Mathematisches Institut, Universit\"at zu K\"oln,
Weyertal 86--90, 50931 K\"oln, F.R.G.}
\email{geiges@math.uni-koeln.de}

\begin{abstract}
A simple proof is given of the following
result first observed by J.~Adachi: embedded circles tangent to the
standard Engel structure on $\R^4$ are classified, up to
isotopy via such embeddings, by their rotation number.
\end{abstract}

\maketitle

\section{Introduction}

The {\it standard Engel structure} $\D$ on $\R^4$ is the maximally
non-integrable $2$-plane distribution defined, in terms
of Cartesian coordinates $(x,y,z,w)$, by the equations
\[ dz-y\, dx=0\;\;\;\mbox{\rm and}\;\;\; dw-z\, dx=0.\]
In other words, it is the tangent $2$-plane field spanned by the
vector fields
\[ e_1:=\partial_x+y\partial_z+z\partial_w\;\;\;\mbox{\rm and}\;\;\;
e_2:=\partial_y.\]

A {\it horizontal loop} is an embedding $\gamma\co S^1\rightarrow\R^4$
everywhere tangent to~$\D$. If we write $\gamma (s)=
(x(s),y(s),z(s),w(s))$, $s\in S^1$, the condition for $\gamma$ to be
horizontal becomes
\[ z'(s)-y(s)x'(s)=0\;\;\;\mbox{\rm and}\;\;\;
   w'(s)-z(s)x'(s)=0\;\;\;\mbox{\rm for all}\;\; s\in S^1.\]
A {\it horizontal isotopy} is an isotopy via horizontal loops.
The {\it rotation number} $\rot (\gamma )$
of a horizontal loop $\gamma$ is the number of
complete turns of the velocity vector $\gamma'(t)\in\D_{\gamma (t)}$,
as we once traverse the loop in positive direction,
relative to the trivialisation of $\D$ given by $e_1,e_2$. The rotation
number is clearly invariant under horizontal isotopies, and it
is easy to show by examples (see Section~\ref{section:proof})
that every integer can be realised as the
rotation number of a horizontal loop.

The following theorem was proved by J. Adachi in~\cite{adac07}:

\begin{thm}
Two horizontal loops are horizontally isotopic if and only if their
rotation numbers agree.
\end{thm}

Adachi proves this theorem by studying the image of $\gamma$ under the
projection
$(x,y,z,w)\mapsto (x,w)$.
(Beware
that I have interchanged $y$ and $w$ in the definition of $\D$ compared
with Adachi's notation. This is more in line with the usual conventions
as regards the contact geometric aspects of our
discussion.)
He determines the
`Reidemeister moves' in this projection and then reduces
the proof to the corresponding classification of topologically
trivial Legendrian knots in standard contact $3$-space, due to
Ya.~Eliashberg and M.~Fraser~\cite{elfr98}.

The purpose of the present note is to show that a much shorter proof
can be given by using the projection $(x,y,z,w)\mapsto (x,z)$
instead.

\section{Horizontal Loops, Legendrian Immersions, and Fronts}
Let $\gamma (s)=(x(s),y(s),z(s),w(s))$, $s\in S^1$, be
a horizontal loop.
Notice that $x'(s)=0$ implies $z'(s)=w'(s)=0$ and hence --- $\gamma$
being an embedding --- $y'(s)\neq 0$.
This means that
\[ \ogamma (s):= (x(s),y(s),z(s)),\;\; s\in S^1,\]
defines a Legendrian immersion into $\R^3$ with its standard
contact structure
\[ \xi := \ker (dz-y\, dx).\]

The {\it rotation number} $\rot (\ogamma )$ --- in the contact
geometric sense --- of such a Legendrian immersion is defined as the
number of complete turns made by $\ogamma'(s)\in\xi_{\ogamma (s)}$
as the loop $\ogamma$ is traversed once in positive direction,
relative to the trivialisation of $\xi$ given by the vector fields
\[ \overline{e}_1:=\partial_x+y\partial_z\;\;\;\mbox{\rm and}
\;\;\; \overline{e}_2:=\partial_y.\]
So it is obvious that $\rot (\gamma )=\rot (\ogamma )$.

The rotation number $\rot (\ogamma )$ is invariant under Legendrian regular
homotopies, i.e.\ $C^1$-homotopies via Legendrian immersions. Moreover,
one can prove by elementary methods that the map $[\ogamma ]\mapsto
\rot (\ogamma )$ defines a one-to-one correspondence between Legendrian
regular homotopy classes of Legendrian immersions $\ogamma\co
S^1\rightarrow (\R^3,\xi )$ on the one hand, and the integers on
the other; see \cite[Thm.~6.3.10]{geig08} or~\cite{geig08a}.

The {\it front projection} $\ogF$ of the Legendrian immersion
$\ogamma$ is the curve
\[ \ogF (s):= (x(s),z(s)), \;\; s\in S^1.\]
Generically, this is an immersed curve with semi-cubical cusps,
but without vertical tangencies~\cite[Section~3.2]{geig08}.
We call a planar curve of this type a {\it front}.

The coordinate $y(s)$ can be recovered as the slope of the front
projection:
\[ y(s)=\frac{z'(s)}{x'(s)}=\frac{dz}{dx}(s).\]
The coordinate $w(s)$ can be recovered as an `area integral':
\[ w(s)-w(s_0)=\int_{s_0}^s z(\sigma )x'(\sigma )\, d\sigma =\int z\, dx.\]
This equation defines the horizontal lift of a Legendrian immersion
$\ogamma$ even when its front projection $\ogamma_F$ is singular, e.g.\ 
during the first Legendrian Reidemeister move~\cite[Figure~7]{geig08a}.

Thus, the condition for an arbitrary front $\ogF$ to lift to a
horizontal {\it immersion}  $\gamma\co S^1\rightarrow (\R^4,\D )$
is that
\[ \oint_{\ogamma_F} z\, dx =0.\]
In order for the lifted curve to be an {\it embedding}, we need to require
in addition that for any two distinct points $s_0,s_1\in S^1$ where
the front has a self-tangency, i.e.\
\[ x(s_0)=x(s_1),\;\;\; z(s_0)=z(s_1),\;\;\;\mbox{\rm and}\;\;\;
\frac{dz}{dx}(s_0)=\frac{dz}{dx}(s_1)\]
--- where, in other words, the lifted Legendrian immersion
$\ogamma$ has a self-intersection ---, we have
\[ \int_{\ogF (s_0)}^{\ogF (s_1)} z\, dx\neq 0.\]
\section{Proof of the Theorem}
\label{section:proof}
\begin{figure}[h]
\centering
\includegraphics[scale=0.46]{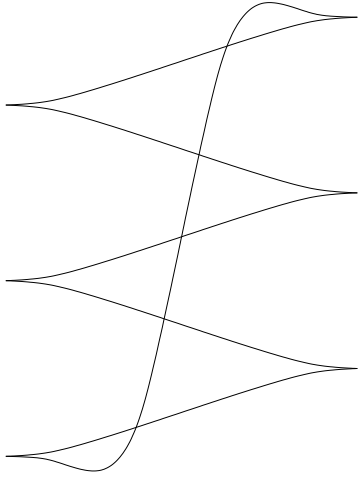}
  \caption{Front of a horizontal loop with $\rot =\pm 3$.}
  \label{figure:model-pm3}
\end{figure}

Figure \ref{figure:model-pm3} shows a front $\ogF$ in the $xz$-plane
whose lift $\ogamma\co S^1\rightarrow (\R^3,\xi )$ is a Legendrian embedding
with $\rot (\ogamma )=\pm 3$, depending on the choice of orientation,
for the rotation number can be computed from the front projection
via the formula
\[ \rot (\ogamma )=\frac{1}{2} (c_- -c_+),\]
with $c_{\pm}$ denoting the number of cusps oriented upwards or
downwards, respectively~\cite[Prop.~3.5.19]{geig08}.

\begin{figure}[h]
\centering
\includegraphics[scale=0.46]{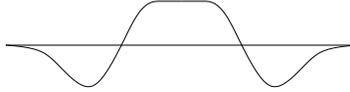}
  \caption{Front of a horizontal loop with $\rot =0$.}
  \label{figure:model-0}
\end{figure}

Furthermore, the area condition discussed in the preceding section
being satisfied, this lifts to a horizontal loop
$\gamma\co S^1\rightarrow (\R^4,\D )$ with $\rot (\gamma )=\pm 3$.
For other non-zero rotation numbers, the picture is entirely analogous;
a front corresponding to a horizontal loop with rotation number zero
is depicted in Figure~\ref{figure:model-0}.

Now consider two horizontal loops $\gamma_0,\gamma_1\co S^1\rightarrow
(\R^4,\D )$ with $\rot (\gamma_0)=\rot (\gamma_1)$. The corresponding
Legendrian immersions $\ogamma_0,\ogamma_1\co S^1\rightarrow (\R^3,\xi )$
satisfy $\rot (\ogamma_0)=\rot (\ogamma_1)$, so there is a Legendrian
regular homotopy between them. As the proof in \cite[p.~312]{geig08}
or \cite[Prop.~4]{geig08a} shows, this can be realised by a homotopy of
fronts $\ogamma_{t,\mathrm{F}}$, $t\in [0,1]$, where at finitely many times
$t_i\in [0,1]$ the front $\ogamma_{t_i,\mathrm{F}}$ has
either a single self-tangency or a singular point (during a first
Legendrian Reidemeister move).
The fronts $\ogamma_{0,\mathrm{F}}$ and 
$\ogamma_{1,\mathrm{F}}$ --- being the projections of horizontal
loops --- satisfy the area conditions, and one can
easily adjust the homotopy of fronts such that all
$\ogamma_{t,\mathrm{F}}$, $t\in [0,1]$, satisfy these conditions,
possibly at the cost of creating further tangencies.
Then the lifted curves $\gamma_t\co S^1\rightarrow (\R^4,\D )$
give the desired horizontal isotopy.

\end{document}